\documentclass[a4paper, 11pt]{amsart}
\usepackage{amsfonts, amsmath, amssymb, amsthm,mathtools, url,dsfont}
\usepackage{graphicx}
\usepackage[english] {babel}
\usepackage{enumitem}
\usepackage[draft]{fixme}
\usepackage{array}
\usepackage{cellspace}
\newtheorem{theorem}{Theorem}

\theoremstyle{definition}

\newcommand{\ve}{\varepsilon}

\usepackage{subcaption}
\usepackage{graphicx}

\begin{document}

\begin{abstract}
In 2004 the second author of the present paper proved that a point set in $[0,1]^d$ which has star-discrepancy at most $\ve$ must necessarily consist of at least $c_{\textup{abs}} d \ve^{-1}$ points. Equivalently, every set of $n$ points in $[0,1]^d$ must have star-discrepancy at least $c_\textup{abs} d n^{-1}$. The original proof of this result uses methods from Vapnik--Chervonenkis theory and from metric entropy theory. In the present paper we give an elementary combinatorial proof for the same result, which is based on identifying a sub-box of $[0,1]^d$ which has $\approx d$ elements of the point set on its boundary. Furthermore, we show that a point set for which no such box exists is rather irregular, and must necessarily have a large star-discrepancy.
\end{abstract}

\title[Irregularities of distributions and extremal sets]{Irregularities of
distributions and extremal sets in combinatorial complexity theory}

\author{Christoph Aistleitner} 
\address{Christoph Aistleitner, Institute of Analysis and Number Theory,
TU Graz, Austria}
\email{aistleitner@math.tugraz.at}

\author{Aicke Hinrichs} 
\address{Aicke Hinrichs, Institute of Analysis, University Linz, Austria}
\email{aicke.hinrichs@jku.at}

\thanks{The first author is supported by the Austrian Science Fund (FWF),
projects F5507-N26, I1751-N26 and Y-901-N35.}

\maketitle

\section{Introduction and statement of results}

Let $\mathcal{A}^*$ denote the class of all axis-parallel boxes in $[0,1]^d$
which have one vertex at the origin. The \emph{star-discrepancy} of a point set
$\mathbf{x}_1, \dots, \mathbf{x}_n \in [0,1]^d$ is defined as 
$$
D_n^*(\mathbf{x}_1, \dots, \mathbf{x}_n) =
\sup_{A \in \mathcal{A}^*} \left| \frac{1}{n} \sum_{k=1}^n \mathds{1}_A(\mathbf{x}_k) -
\textup{vol} (A) \right|,
$$
where $\mathds{1}_A$ denotes the indicator function of $A$.\footnote{It does not make any difference, but for convenience we will assume in this paper that the boxes in $\mathcal{A}^*$ are closed. We will also allow point sets to contain identical points, so strictly speaking our point sets are not sets, but multi-sets.} The notion of the star-discrepancy is crucial for the \emph{Quasi-Monte Carlo
integration} method, in which the integral over $[0,1]^d$ of a $d$-variate function $f$ is
approximated by the average $\frac{1}{n} \sum_{k=1}^n f(\mathbf{x}_k)$.
Famously, by the Koksma--Hlawka inequality the error in this numerical
integration method can be estimated by the product of the variation of $f$ in an
appropriate sense and the star-discrepancy of the set of sampling points
$\mathbf{x}_1, \dots, \mathbf{x}_n$. More information on this topic can be found
in the classical monographs \cite{dpd,dts,knu}.\\

The most famous open problem in discrepancy theory concerns the
necessary degree of irregularity of a point distribution in the multidimensional
unit cube. More precisely, the problem asks for the smallest possible
order of the discrepancy of a set of $n$ points in $[0,1]^d$, which was partially answered by the
celebrated results of Roth~\cite{roth} and Bilyk--Lacey--Vaghar\-shakyan
\cite{blv} (see \cite{bilyk} for a survey) on the one hand and by many
constructions of so-called low-discrepancy point sets (see \cite{dpd}) on the
other hand. In the formulation of this problem it is understood that $d$ is
fixed and $n \to \infty$. Another important open problem, which recently has received some attention, asks for the
order of the \emph{inverse of the star-discrepancy}: given a number $\varepsilon>0$, what is
the minimal cardinality $n^*(\varepsilon,d)$ of a point set in $[0,1]^d$ achieving
star-discrepancy at most $\varepsilon$? This problem can also been seen as an
irregularities-of-distributions problem, but one where the role of the
\emph{simultaneous} dependence of the 
minimal size of the discrepancy on both $d$ and $n$ is emphasized.\\

Concerning the inverse of the discrepancy, it is known that
\begin{equation} \label{inv}
n^*(\varepsilon,d) \leq c_{\textup{abs}} d \varepsilon^{-2} 
\end{equation}
from a fundamental paper of Heinrich--Novak--Wasil\-kowski--Wo{\'z}nia\-kowski
\cite{hnww}, and that
\begin{equation} \label{hin}
n^*(\varepsilon,d) \geq c_{\textup{abs}} d \varepsilon^{-1}  \qquad \qquad (\varepsilon < \varepsilon_0)
\end{equation}
due to a result of the second author of the present paper
\cite{hinrichs}.\footnote{Throughout this paper, $c_{\textup{abs}}$ denotes
positive absolute constants, not always the same.} Thus the inverse of the
discrepancy depends \emph{linearly} on the dimension $d$, while the dependence
on $\varepsilon$ constitutes an important open problem. Novak and
Wo{\'z}nia\-kowski conjectured that the exponent 2 of $\varepsilon^{-1}$ in \eqref{inv} is
optimal. In~\cite[p.~63]{nw2} they write:
\begin{quote}
How about the dependence on $\varepsilon^{-1}$? This is open and seems to be a difficult
problem. [\dots] We think that as long as we consider upper bounds of the form
$n^*(\varepsilon,d) \leq c_{\textup{abs}} d^k \varepsilon^{-\alpha}$, the exponent $\alpha \geq
2$ and 2 cannot be improved. 
\end{quote}
See also~\cite[Open problem 7]{nw1} and~\cite[Problem 3]{hein}.\\

Note that \eqref{inv} and \eqref{hin} can be formulated in a different, alternative form.
Equation \eqref{inv} is equivalent to saying that for all $d$ and $n$ there
exist $\mathbf{x}_1, \dots, \mathbf{x}_n \in [0,1]^d$ such that
$$
D_n^* (\mathbf{x}_1, \dots, \mathbf{x}_n ) \leq c_{\textup{abs}}
\frac{\sqrt{d}}{\sqrt{n}},
$$
while \eqref{hin} is equivalent to the statement that for all $\mathbf{x}_1,
\dots, \mathbf{x}_n \in [0,1]^d$ we have
\begin{equation} \label{hin2}
D_n^*(\mathbf{x}_1, \dots, \mathbf{x}_n) \geq c_{\textup{abs}} \frac{d}{n}
\qquad \qquad (n \geq c_{\textup{abs}} d).
\end{equation}

The proof of \eqref{hin} in \cite{hinrichs} uses methods from combinatorial
complexity theory (more precisely, Vapnik--Chervonenkis theory) together with methods from 
metric entropy theory. The purpose of the present paper is twofold. On the one
hand, we want to give an elementary proof of \eqref{hin}, in the spirit of the
``cheap proof'' which will be sketched below. On the other hand, we will use
Vapnik--Chervonenkis theory (VC theory) and metric entropy theory in order to
show that point sets which prohibit an application of the ``cheap proof'' must
necessarily have a rather simple combinatorial structure from the point of view
of VC theory, and must consequently have particularly large discrepancy.\\

The idea of the ``cheap proof'' is very simple. Let $\mathbf{x}_1, \dots,
\mathbf{x}_n$ be points in $[0,1]^d$. Find a box $A \in \mathcal{A}^*$ such that
$d$ points of $\mathbf{x}_1, \dots, \mathbf{x}_n$ are situated on the ``right
upper'' boundary of $A$ (that is, on one of the faces which are \emph{not}
adjacent to the origin). Let $A_1$ be a box which is just a little bit
smaller than $A$ and let $A_2$ a box which is just a little bit larger than $A$.
Then the volumes of $A_1$ and $A_2$ are essentially equal, while the difference
in points is at least $d$. Thus the star-discrepancy of $\mathbf{x}_1, \dots,
\mathbf{x}_n$ is at least $d/n$.\\

The problem with the ``cheap argument'' clearly is that it is not always
possible to find a box $A$ which has $d$ points on its boundary -- see for
example the point set on the left-hand side of Figure 1 below. However, our proof of Theorem
\ref{th1} shows that a slight modification of the ``cheap argument'' can
actually be successfully implemented. Furthermore, as Theorem \ref{th2} will
show, a point set which does not allow the ``cheap argument'' must have a very
strong internal structure, and in particular must have a small combinatorial
complexity in the sense of VC theory.

\begin{theorem} \label{th1}
Let $\mathbf{x}_1, \dots, \mathbf{x}_n$ be points in $[0,1]^d$. Then
$$
D_n^*(\mathbf{x}_1, \dots, \mathbf{x}_n) \geq \frac{d}{12n},
$$
provided that $n \geq 250 d$.
\end{theorem}

From Theorem \ref{th1} we can deduce that
$$
n^*(\varepsilon,d) \geq \frac{d \varepsilon^{-1}}{12}, \qquad \qquad \left(\varepsilon < \frac{1}{3000}
\right).
$$
The constants appearing in Theorem \ref{th1} may be compared with those given in
\cite{hinrichs} for \eqref{hin}, where it was shown that
$$
D_n^*(\mathbf{x}_1, \dots, \mathbf{x}_N) \geq \frac{d}{32 e^2 n}, \qquad \qquad (n \geq d),
$$
with $32 e^2 \approx 236$. However, the reason for writing the present paper was
not to improve the numerical constants in \eqref{hin}; rather, the purpose of
this paper is to share some observations which we consider interesting.\\

The following theorem states, informally speaking, that a point configuration
which does not allow one to apply the ``cheap argument'' must necessarily have a
small combinatorial complexity, and consequently must have a large discrepancy.
In other words, either the cheap argument is applicable straightforward or the
point configuration must have even larger discrepancy especially because it
prohibits the application of the cheap argument. In the statement of the theorem, as the ``right
upper'' boundary of an anchored axis-parallel box $A = [\mathbf{0},\mathbf{a}]$ we mean the union of all those $(d-1)$-dimensional faces of $A$ which are adjacent to the (``right upper'') point $\mathbf{a}$.

\begin{theorem} \label{th2}
Let $\mathbf{x}_1, \dots, \mathbf{x}_n$ be points in $[0,1]^d$, and assume that
it is not possible to find a box in $\mathcal{A}^*$ such that the right upper
boundary of this box contains at least $d/4$ of these points. Assume also that $n \geq d$. Then
$$
D_n^*(\mathbf{x}_1, \dots, \mathbf{x}_n) \geq \frac{d^{3/4}}{372 n^{3/4}}.
$$
\end{theorem}

We finish the introduction with a discussion on the applicability of the ``cheap
proof'' and on the combinatorial complexity of point sets. Here combinatorial
complexity refers to the cardinality of the set
\begin{equation} \label{card}
\bigl\{ A \cap \{ \mathbf{x}_1, \dots, \mathbf{x}_n \}:~A \in \mathcal{A}^*
\bigr\}
\end{equation}
(this is a set of subsets of $\{\mathbf{x}_1, \dots, \mathbf{x}_n \}$). Since
the class of anchored axis-parallel boxes $\mathcal{A}^*$ is a
\emph{Vapnik--Chervonenkis class} (VC class) of index $d$, the cardinality of the set
\eqref{card} can be bounded by the Sauer--Shelah lemma, which asserts that this
cardinality is at most
\begin{equation} \label{sauer}
\sum_{i=0}^d \binom{n}{i};
\end{equation}
this is one of the main ingredients of the ``entropy argument'' in the proof of
\eqref{hin} in \cite{hinrichs} (for definitions in the context of VC theory, see Section \ref{complex}). This entropy argument gives better (that is,
larger) lower bounds for the discrepancy the smaller the
cardinality of \eqref{card} can be shown to be. As we will show in Section
\ref{complex} below, a point set which prohibits the application of the ``cheap
argument'' by having the property stated in the assumption of Theorem \ref{th2}
must have a very small combinatorial complexity in the sense that the
cardinality of \eqref{card} is much smaller than what could be deduced from the
Sauer--Shelah lemma. However, in turn, if an improvement of the Sauer--Shelah lemma is
not possible since the point set does not satisfy the assumptions of Theorem
\ref{th2}, then obviously the ``cheap argument'' is applicable to this point
set. So two competing forces are at work here, both of which lead to a large
discrepancy in one way or the other.\\

To illustrate the situation we present two extremal point sets in the picture
below (unfortunately the picture is restricted to the less instructive
two-dimensional case). The point set (A) on the left-hand side is extremal in
the sense that it prohibits the application of the ``cheap argument'' -- there
is no anchored axis-parallel box which has two elements of the point set on its
right-upper boundary. On the other hand, the point set (A) has very low
complexity in the sense of VC theory: there are 9 points, and the cardinality of
\eqref{card} is obviously 10 (note that the empty set also counts), which is
smallest possible (unless we allow points to coincide). In contrast, for the point set (B) the cardinality of
\eqref{card} can be calculated to be 46, which is $\binom{9}{0} + \binom{9}{1} +
\binom{9}{2}$ and thus by \eqref{sauer} is largest possible. On the other hand,
the ``cheap argument'' is obviously applicable to this point set, and actually
there is a very large number of boxes which have two elements of (B) on their
``right upper'' boundary.

\begin{figure}
\centering
\begin{subfigure}{.5\textwidth}
  \centering
  \includegraphics[width=0.85 \linewidth]{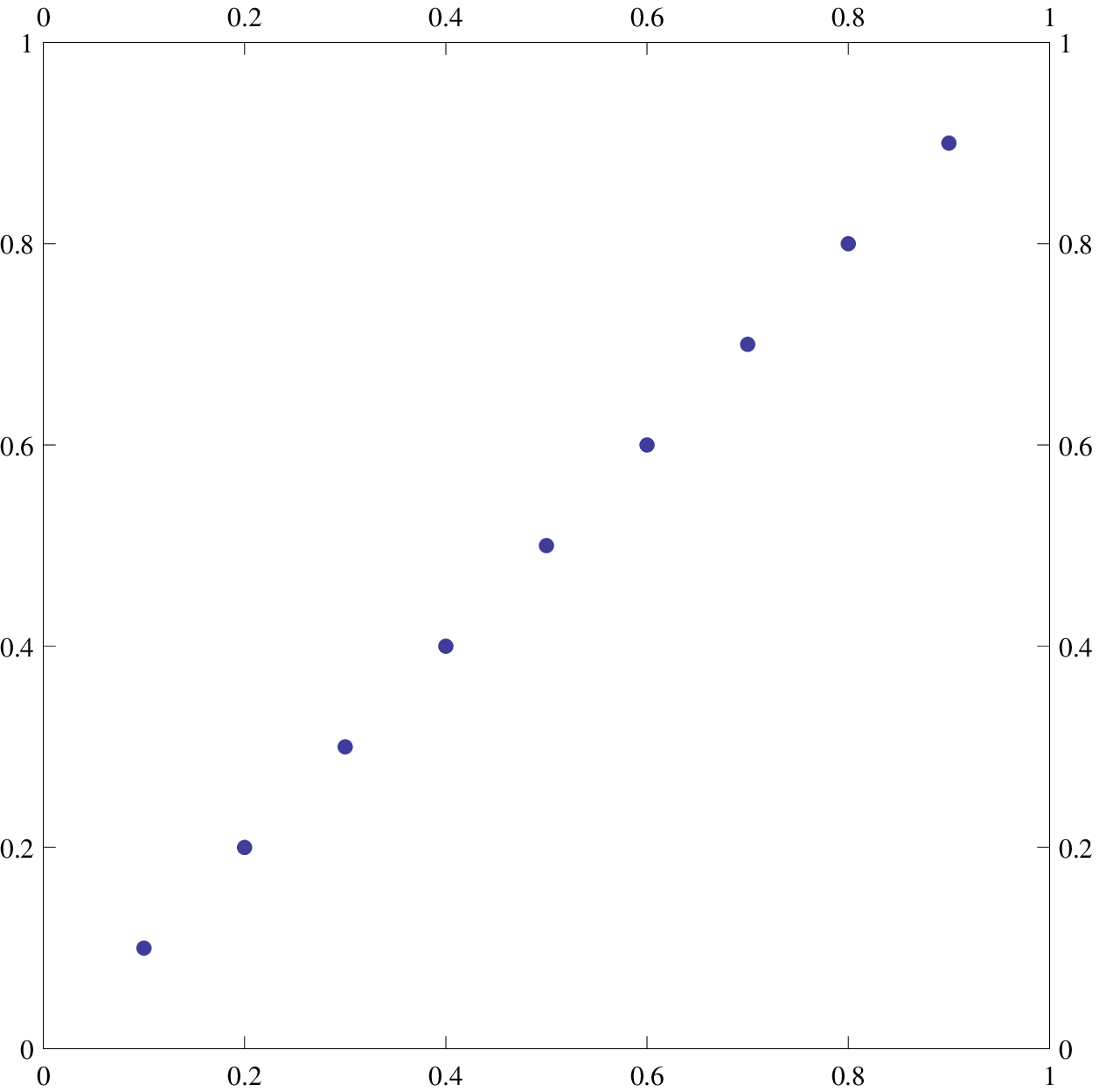}
  \caption{}
  \label{fig:sub1}
\end{subfigure}%
\begin{subfigure}{.5\textwidth}
  \centering
  \includegraphics[width=0.85 \linewidth]{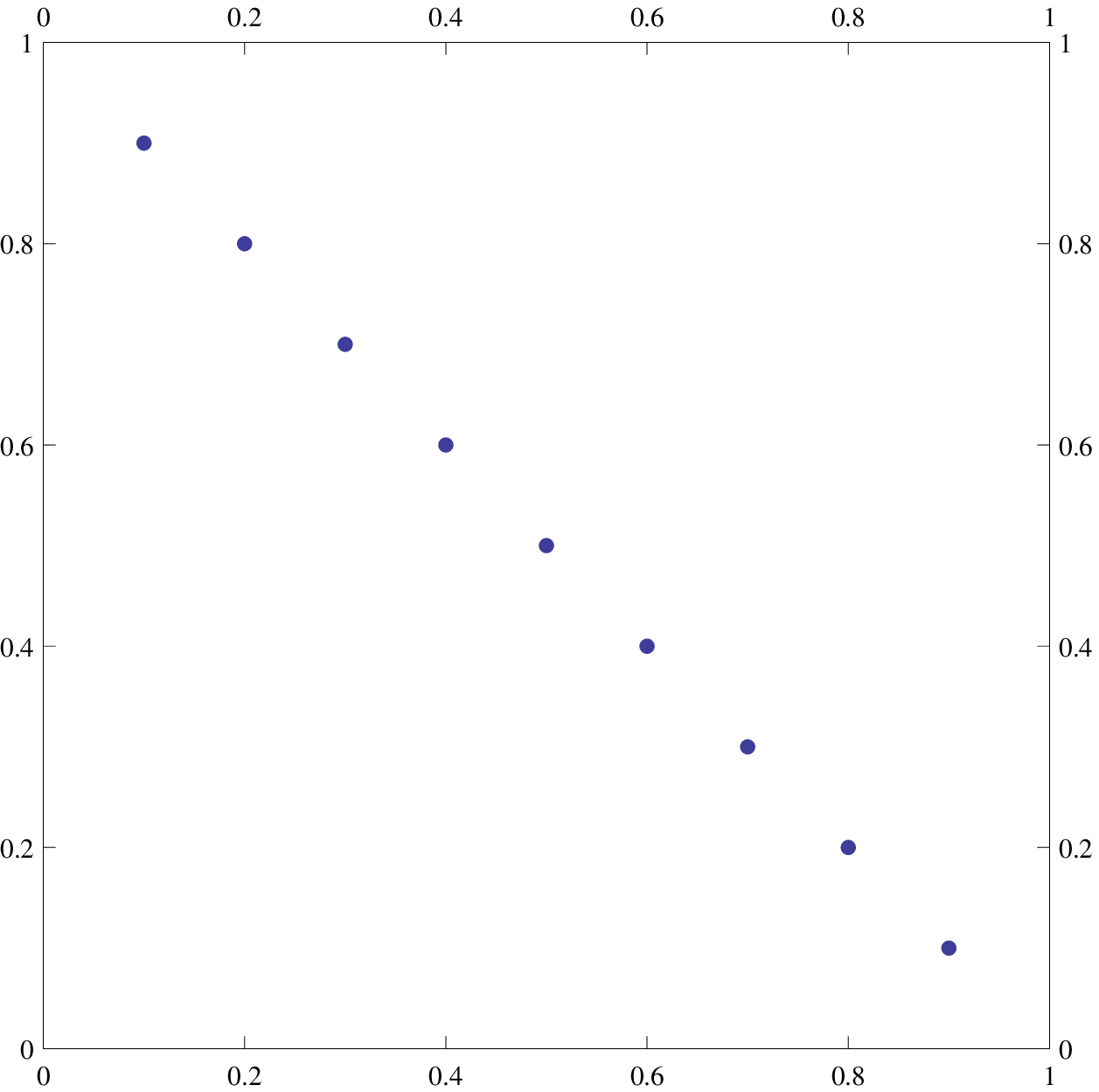}
  \caption{}
  \label{fig:sub2}
\end{subfigure}
\caption{Two extremal point sets. Note that the point set (A) is exactly the ``opposite'' of point set (B).}
\label{fig:test}
\end{figure}

The outline of the remaining part of this paper is as follows. In Section 
\ref{cheap} we use the ``cheap argument'' to prove Theorem \ref{th1} in the case
$n \geq 2 e d^2$, which is simpler than the general case and particularly
instructive. In Section \ref{cheapg} we use the ``cheap argument'' to prove
Theorem \ref{th1} in the general case. In Section \ref{complex} we introduce the
necessary notions from VC theory and prove Theorem
\ref{th2}.

\section{The cheap proof of Theorem \ref{th1} in the case $n \geq 2 e d^2$}
\label{cheap}

In this section we will prove Theorem \ref{th1} under the additional assumption that
$n \geq 2 e d^2$, since in this case the proof is particularly simple. It also
illustrates the idea of the proof in the general case, which, however, requires
a more careful reasoning.

We may assume that $d \geq 2$. Let $\mathcal{P}=\{\mathbf{x}_1, \dots, \mathbf{x}_n\}$ be a point set in $[0,1]^d$
and let $\kappa=1-\frac1d$. Now consider the boxes $A,B \in \mathcal{A}^*$ given by
$$ 
A=[0,1] \times [0,\kappa]^{d-1} \qquad \text{and} \qquad B = [0,\kappa]^d.
$$ 
Observe that $B \subset A$ and that
\begin{equation}\label{ij}
 \textup{vol~} (A \setminus B) = \frac{1}{d} \left( 1-\frac{1}{d} \right)^{d-1} > \frac{1}{ed} \geq \frac{2d}{n},
\end{equation}
where we used the assumption that $n \geq 2 e d^2$.
If $A \cap \mathcal{P} = B \cap \mathcal{P}$ then we find
\begin{align*}
 2 D_n^* (\mathbf{x}_1, \dots, \mathbf{x}_n) & \geq \left( \textup{vol~} (A) - \frac{\# (A \cap \mathcal{P})}{n}\right) - \left( \frac{\# (A \cap \mathcal{P})}{n} - \textup{vol~} (B) \right) \\
  &= \textup{vol~} (A \setminus B)  \geq \frac{2d}{n},
\end{align*}
which implies
\begin{equation} \label{eq:dd}
 D_n^* (\mathbf{x}_1, \dots, \mathbf{x}_n) \ge \frac{d}{n}.
\end{equation}
On the other hand, if $A \cap \mathcal{P} \neq B \cap \mathcal{P}$, then there exists a point $\mathbf{y}_1 = (y_1^{(1)},\dots,y_1^{(d)}) \in \mathcal{P} $ in $A \setminus B$, i.e.
$$
 y_1^{(1)} > \kappa \qquad \text{and} \qquad  y_1^{(k)} \le \kappa \ \text{ for } k \neq 1. 
$$
Arguing similarly for the other coordinates, we either have already proved \eqref{eq:dd} or we find $y_1,\dots,y_d\in \mathcal{P} $ such that 
$$
 y_j^{(j)} > \kappa \qquad \text{and} \qquad  y_j^{(k)} \le \kappa \ \text{ for } k \neq j. 
$$
Obviously, all these points are distinct and contained in the right upper boundary of the box
$$
\left[0, y_1^{(1)} \right] \times \left[0, y_2^{(2)} \right] \times \dots \times
\left[0, y_d^{(d)} \right].
$$
Thus we have found a box which contains at least $d$ elements of $\mathcal{P}$ on its right upper boundary, and the ``cheap argument'' from above shows that \eqref{eq:dd} holds in this case.
This proves Theorem \ref{th1} (with the value
$1$ instead of $1/12$ for the constant) in the case $n \geq 2 e d^2$.

\section{The ``cheap proof'' of Theorem \ref{th1} in the general case}
\label{cheapg}

The proof of Theorem \ref{th1} in the general case uses the same idea as the
proof in
the previous section; however, it requires a slightly more complicated combinatorial
argument. The reason why the argument from the previous section fails is that
the last inequality of \eqref{ij} is no longer true, so that we can no longer
guarantee that every box of the type $A \setminus B$ contains a point. Consequently we will use a slightly different construction, and distinguish between several cases.\\

As the reader will see our proof contains several numerical parameters, such as
the number 25 in the definition of $\kappa$ below. Of course we have chosen
parameters which give a reasonable result. However, optimizing these parameter
is fairly complicated, and we do not claim that we have found the optimal ones.
Also, with this method there is a trade-off between the two constants appearing in
the statement of the theorem, which are 1/12 and 250 in our formulation of Theorem \ref{th1}.
Decreasing one of them would possibly increase the other, and vice versa.
In particular, we checked that the theorem also holds with constants 1/20 and 40.\\

We fix the point set $\mathcal{P} = \{ \mathbf{x}_1, \dots, \mathbf{x}_n \}$ and abbreviate 
$$
 D_n^* = D_n^* (\mathbf{x}_1, \dots, \mathbf{x}_n)
$$
and $[d]=\{1,\dots,d\}$.
The trivial bound $D_n^* \ge \frac{1}{2n}$ from the one-dimensional case already proves the theorem in the case $d\le 6$. So we may and do assume that $d \ge 7$.\\

We will need the reverse Bernoulli-type inequality
\begin{equation}\label{eq:Bern}
 (1-x)^q \ge 1- \frac{21}{20} q x \qquad \text{for } 0\le x\le \frac{1}{10} \text{ and } \frac{1}{7} \le q \le \frac{1}{4}.
\end{equation}
This inequality can be easily checked numerically and, what is more tedious, can be proved by elementary analysis.\\

With
$$
\kappa = \left(1-\frac{25d}{n}\right)^{1/d}
$$ 
we partition the point set $\mathcal{P}$ into subsets according to how many coordinates of the considered point are at least $\kappa$:
\begin{align*}
 \mathcal{P}_0 &= \big\{ \mathbf{x} \in \mathcal{P} \,:\, x^{(j)} \le \kappa \text{ for all } j \in [d] \big\}, \\
 \mathcal{P}_1 &= \big\{ \mathbf{x} \in \mathcal{P} \,:\, x^{(j)} > \kappa \text{ for exactly one } j \in [d] \big\},\\ 
 \mathcal{P}_2 &= \big\{ \mathbf{x} \in \mathcal{P} \,:\, x^{(j)} > \kappa \text{ for at least two } j \in [d] \big\}.
\end{align*}
Furthermore, let
$$
 \mathcal{C} = \big\{ j \in [d] \,:\, x^{(j)} > \kappa \text{ for some } \mathbf{x} \in   \mathcal{P}_1 \big\}
$$
be the set of coordinates where at least one point in $\mathcal{P}_1 $ has its largest coordinate.

We now distinguish between three cases.\\

\emph{Case 1.} Assume that $\# \mathcal{C} \ge \frac{d}{6}$. \\

This is the simple case, when the ``cheap proof'' is directly applicable. Since every $\mathbf{x}\in \mathcal{P}_1$ has exactly one coordinate $j$ with $x^{(j)}> \kappa$, there exists a box $A \in \mathcal{A}^*$ that contains $\# \mathcal{C} \ge \frac{d}{6}$ points of $\mathcal{P}_1 \subseteq \mathcal{P}$ on its right upper boundary. Hence the ``cheap proof'' shows that 
$$
 D_n^* \geq \frac{d}{12n}
$$
in this case.\\

\emph{Case 2.} Assume that $\# \mathcal{C} < \frac{d}{6}$ and $\# \mathcal{P}_1 \geq \frac{107d}{24}$.\\

In this case there exist many points having exactly one large coordinate, but the ``cheap proof'' is not applicable since many of these points share the same few coordinate indices where they have their large coordinate. However, since too many points are located close to just a few right upper faces of the unit cube, there must also exist a large sub-box of $[0,1]^d$ (avoiding the proximity of these faces) which does not contain enough points.\\

More precisely, the box
$$
 A = [0,\kappa]^{\mathcal{C}} \times [0,1]^{[d] \setminus \mathcal{C}}.
$$
has volume
$$
 \textup{vol} (A) = \kappa^{\# \mathcal{C}} \ge \kappa^{d/6} = \left(1-\frac{25d}{n}\right)^{1/6} \ge 1- \frac{35d}{8n},
$$
where the last inequality follows from \eqref{eq:Bern} and the assumption $n\ge 250 d$. By definition of the sets  $\mathcal{P}_i$, we have
$$
 \# (A \cap \mathcal{P}) \le \# \mathcal{P}_0 +  \# \mathcal{P}_2  = n - \# \mathcal{P}_1 \le n-\frac{107d}{24}.
$$ 
This implies
$$
 D_n^* \ge \textup{vol} (A) - \frac{\# (A \cap \mathcal{P})}{n} \ge \left(\frac{107}{24} - \frac{35}{8}\right) \frac{d}{n} = \frac{d}{12n} %> \frac{d}{10n}
$$
also in this case.\\

\emph{Case 3.} Assume that $\# \mathcal{P}_1 < \frac{107d}{24}$.\\

This is the most tricky case. Since the cardinality of $\mathcal{P}_1$ is small, the cardinality of $\mathcal{P}_2$ must be large. Thus we have a relatively large number of points which have multiple large coordinates, which means that these points cannot be assigned to different faces of the unit cube (as in Section \ref{cheap} or as in Case 1) but that they are rather located in ``corners'' of the unit cube. Thus the ``cheap proof'' is not applicable. However, since many points are located in corners, this means that we can identify a large sub-box of $[0,1]^d$, reaching all the way from 0 to 1 in many coordinates, which avoids these corners and contains an insufficient number of points of $\mathcal{P}$.\\

To give a detailed proof in this case, first we consider  $A = [0,\kappa]^d$ which has volume $\kappa^d = 1-\frac{25d}{n}$ and contains exactly those points of $\mathcal{P}$ that are in  $\mathcal{P}_0$, that is
$$
 \# (A \cap \mathcal{P}) = \# \mathcal{P}_0 = n - \# \mathcal{P}_1 - \# \mathcal{P}_2.
$$
Then it follows from 
$$
  1- \frac{\# \mathcal{P}_1 + \# \mathcal{P}_2}{n}  - \textup{vol} (A) =\frac{ \# (A \cap \mathcal{P})}{n} - \textup{vol} (A) \le D_n^*
$$
and from the assumption $\# \mathcal{P}_1 < \frac{107d}{24}$ that
\begin{equation} \label{eq:M}
 M:= \# \mathcal{P}_2 \geq \frac{493d}{24} - n D_n^*.
\end{equation}

We now set up an inductive procedure to produce a large box which contains few points by successively removing points of $\mathcal{P}_2$. Let
$$
  S_0=\mathcal{P}_2, \ R_0 = \emptyset, \ \mathcal{C}_0 = \emptyset, \ m_0=\# R_0=0.
$$ 
Now assume that $S_{k-1},R_{k-1} \subset \mathcal{P}_2, \mathcal{C}_{k-1} =\{j_1,\dots,j_{k-1}\}\subset [d]$ and $ m_{k-1}= \# R_{k-1}$ are already defined. By definition of $\mathcal{P}_2$ and double counting we have
\begin{align*}
 \sum_{j\in [d] \setminus \mathcal{C}_{k-1}} \# \big\{ \mathbf{x} \in S_{k-1} \,:\, x^{(j)} > \kappa \big\} 
 &=
 \sum_{ \mathbf{x} \in S_{k-1}} \# \big\{  j\in [d] \setminus \mathcal{C}_{k-1}\,:\, x^{(j)} > \kappa \big\} \\
 &\geq
 2 \#  S_{k-1}.
\end{align*}
Therefore, as long as $S_{k-1} \neq \emptyset$, we find $j_k \in [d] \setminus \mathcal{C}_{k-1}$ such that
$$
 R_k =  \big\{ \mathbf{x} \in S_{k-1} \,:\, x^{(j_k)} > \kappa \big\}
$$
satisfies
$$
 m_k = \# R_k \ge \frac{2 \#  S_{k-1}}{\# [d] \setminus \mathcal{C}_{k-1}}
 \geq \frac{2 \#  S_{k-1}}{d} = \frac{2}{d} \left[ M - \sum_{h=1}^{k-1} m_h\right].
$$
To complete the inductive construction, let
$$
 \mathcal{S}_{k} = \mathcal{S}_{k-1} \setminus R_k \quad \text{and} \quad  
 \mathcal{C}_{k} = \mathcal{C}_{k-1} \cup \{j_k\}.
$$
If $S_{k} = \emptyset$ for some $k<d$, we take $S_h=R_h=\emptyset$ for $h\ge k$ and choose $j_h$ arbitrary among the remaining coordinates. Then the inductive process is defined for $k=0,1,\dots,d$ and all the above estimates hold. Fix $k$ and let $q=\frac{k}{d}$. For the total number of points removed up to step $k$, we then have
$$
 \sum_{h=1}^{k} m_h \ge \frac{2k}{d} \left[ M - \sum_{h=1}^{k} m_h\right] = 2q \left[ M - \sum_{h=1}^{k} m_h\right],
$$
which, by \eqref{eq:M}, implies
$$
 \sum_{h=1}^{k} m_h \ge \frac{2qM}{1+2q} \ge \frac{2q}{1+2q} \left(\frac{493d}{24} - n D_n^*\right).
$$

We now consider the box 
$$
 A = [0,\kappa]^{\mathcal{C}_k} \times [0,1]^{[d] \setminus \mathcal{C}_k}.
$$
which has volume
$$
 \textup{vol} (A) = \kappa^{\# \mathcal{C}_k} \ge \kappa^{k} = \left(1-\frac{25d}{n}\right)^{q} \ge 1- \frac{105q}{4} \cdot \frac{d}{n},
$$
where the last inequality follows from \eqref{eq:Bern} and the assumption $n\ge 250 d$, provided that
$$
 \frac{1}{7} \le q=\frac{k}{d} \le \frac{1}{4}.
$$
Since we assumed $d\geq 7$ and since $d$ is an integer, we can satisfy this condition with the choice $k = \lceil \frac{d}{7} \rceil$.\\

By construction, none of the removed points in $R_1,\dots,R_k$ is contained in $A$, which implies
$$
 \# (A \cap \mathcal{P}) \le n - \sum_{h=1}^{k} m_h \le n - \frac{2q}{1+2q} \left(\frac{493d}{24} - n D_n^*\right).
$$
Hence
$$
 D_n^* \geq  \textup{vol} (A) - \frac{\# (A \cap \mathcal{P})}{n} \geq
 \left( \frac{493 q}{12(1+2q)} - \frac{105q}{4} \right) \frac{d}{n} - \frac{2 q }{1+2q} D_n^*,
$$
which in turn gives
$$
 D_n^* \geq  \frac{1+2q}{1+4q}\left( \frac{493 q}{12(1+2q)} - \frac{105q}{4} \right) \frac{d}{n}
 = \frac{q (89-315q)}{6(1+4q)} \, \frac{d}{n}.
$$
It is easily verified that
$$
 \frac{q (89-315q)}{6(1+4q)} \ge \frac{1}{12}
$$
for $\frac{1}{7} \le q \le \frac{1}{4}$, so that the theorem is also proved in this case.

\section{Point sets which prohibit the ``cheap argument'', and combinatorial
complexity theory} \label{complex}

In Vapnik--Chervonenkis theory (VC theory) the notion of \emph{shattering} plays
a crucial role. Let $S = \{x_1, \dots, x_n\}$ be elements of some set $X$, and
let $\mathcal{C}$ denote a collection of subsets of $X$. We say that
$\mathcal{C}$ \emph{shatters} $S$ if 
$$
\# \left\{A \cap S:~A \in \mathcal{C} \right\} = 2^n;
$$ that is, if using the sets in $\mathcal{C}$ it is possible to
pick out every possible subset from $S$. The \emph{VC index} (or \emph{VC
dimension}) of $\mathcal{C}$ is the smallest $n$ for which no set (of elements of $X$) of cardinality
$n$ exists which is shattered by $\mathcal{C}$. In our setting we have $X =
[0,1]^d$ and $\mathcal{C} = \mathcal{A}^*$, and the VC dimension of
$\mathcal{A}^*$ is $d$.\\

Assume that $\mathcal{C}$ has VC dimension $d$, and that $\# S = n$. Then the Sauer--Shelah lemma
asserts that 
\begin{equation} \label{sauer2}
 \# \left\{A \cap S:~A \in \mathcal{C} \right\} \leq \sum_{i=0}^d \binom{n}{i}
\end{equation}
(and this upper bound is in general optimal). We will sketch a proof of this
lemma in the setting $X = [0,1]^d$ and $\mathcal{C} = \mathcal{A}^*$ below,
since we will use this proof as a blueprint for the key inequality in our proof
of Theorem \ref{th2}.\\

Set 
\begin{equation} \label{ndn}
N(n,d) = \max_{\mathbf{x}_1, \dots, \mathbf{x}_n \in [0,1]^d} \# \left\{ A \cap \{\mathbf{x}_1, \dots, \mathbf{x}_n\}:~A \in \mathcal{A}^*\right\}.
\end{equation}
Let $\mathbf{x}_1, \dots, \mathbf{x}_n$ be any points in $[0,1]^d$. Assume,
without loss of generality, that $\mathbf{x}_1$ has the largest first coordinate
among all these points. Then for a given box $A$ for the intersection $A \cap \{\mathbf{x}_1, \dots, \mathbf{x}_n\}$
there are two possibilities. Either $\mathbf{x}_1 \not\in A$, which means that we ``lose'' one point. Or $\mathbf{x}_1 \in A$, which
means that the first coordinate of the right upper corner of $A$ is at least as large as the first coordinates of all other points as well, and we ``lose'' one dimension as well as the point $\mathbf{x}_1$ (which by construction is always contained in $A$ in this case). Thus
\begin{equation} \label{rec}
N(n,d) \leq N(n-1,d) + N(n-1,d-1), \qquad d,n \geq 2.
\end{equation}
Together with the trivial initial values $N(1,d) = 2$ and $N(n,1) = n+1$ this
leads to a recursion, whose solution gives \eqref{sauer2}. A detailed version of this proof can be found, for example, on page 46 of \cite{mohri}.\\

Now we are ready to prove Theorem \ref{th2}. Let $d$ be fixed. To avoid ambiguities, we write $\mathcal{A}^*(d)$ for the collection of axis-parallel boxes having one vertex at the origin, which are contained in $[0,1]^d$. For points $\mathbf{y}_1, \dots, \mathbf{y}_m$ in $[0,1]^s$, we say that this collection of points has property $\mathbf{P}(r)$ if it is not possible to find a box in $\mathcal{A}^*(s)$ such that the right upper boundary of this box contains at least $r$ of these points (where the term ``right upper boundary'' is defined as in the paragraph before the statement of Theorem \ref{th2}). For the assumptions of Theorem \ref{th2} this means that we start with a set $\mathbf{x}_1, \dots, \mathbf{x}_n$ in $[0,1]^d$ having property $\mathbf{P}(d/4)$.\\

Set
$$
\hat{N} (m,s,r) = \max_{\substack{\mathbf{y}_1, \dots, \mathbf{y}_m \in [0,1]^s,\\ \substack{\mathbf{y}_1, \dots, \mathbf{y}_m \text{ has property $\mathbf{P}(r)$}}}} \# \left\{ A \cap \{\mathbf{y}_1, \dots, \mathbf{y}_m\}:~A \in \mathcal{A}^*(s) \right\}.
$$
Assume that $r < s$, and let $\mathbf{y}_1, \dots, \mathbf{y}_m \in [0,1]^s$ be points having property $\mathbf{P}(r)$. Upon a little reflection this implies that there must exist a point in $\mathbf{y}_1, \dots, \mathbf{y}_m$ which has at least 2 maximal coordinates (that is, coordinate entries which are at least as large as the corresponding coordinate entries of all the other points). Without loss of generality, assume that this point is $\mathbf{y}_m$, and that its coordinates at positions $s-1$ and $s$ are maximal in this sense. Then
\begin{eqnarray}
& & \# \left\{ A \cap \{\mathbf{y}_1, \dots, \mathbf{y}_m\}:~A \in \mathcal{A}^*(s) \right\} \nonumber\\
& = & \# \left\{ A \cap \{\mathbf{y}_1, \dots, \mathbf{y}_{m}\}:~A \in \mathcal{A}^*(s),~\mathbf{y}_m \not\in A \right\} \nonumber\\
& & \quad + \quad \# \left\{ A \cap \{\mathbf{y}_1, \dots, \mathbf{y}_m\}:~A \in \mathcal{A}^*(s),~\mathbf{y}_m \in A \right\} \nonumber\\
& = & \# \left\{ A \cap \{\mathbf{y}_1, \dots, \mathbf{y}_{m-1}\}:~A \in \mathcal{A}^*(s) \right\} \label{term1}\\
& & \quad + \quad \# \left\{ A \cap \{\mathbf{y}_1, \dots, \mathbf{y}_{m-1}\}:~A \in \mathcal{A}^*(s),~\mathbf{y}_m \in A \right\}. \label{term2}
\end{eqnarray}
The term in line \eqref{term1} is clearly dominated by $\hat{N} (m-1,s,r)$. To understand the term in line \eqref{term2}, let $\mathbf{y}^{(s-2)}$ denote the restriction (projection) of a point $\mathbf{y} \in [0,1]^s$ to its first $s-2$ coordinates, and define $A^{(s-2)}$ similarly as a projection of $A$. By construction $\mathbf{y}_m \in A$ implies that the coordinates at positions $s-1$ and $s$ of all the points $\mathbf{y}_1, \dots, \mathbf{y}_{m-1}$ cannot exceed those of the right upper corner of $A$. Thus
\begin{eqnarray}
& & \# \left\{ A \cap \{\mathbf{y}_1, \dots, \mathbf{y}_{m-1}\}:~A \in \mathcal{A}^*(s),~\mathbf{y}_m \in A \right\} \nonumber\\
& = & \# \left\{ A^{(s-2)} \cap \{\mathbf{y}_1^{(s-2)}, \dots, \mathbf{y}_{m-1}^{(s-2)}\}:~A \in \mathcal{A}^*(s),~\mathbf{y}_m \in A \right\} \nonumber\\
& \leq & \# \left\{ A \cap \{\mathbf{y}_1^{(s-2)}, \dots, \mathbf{y}_{m-1}^{(s-2)}\}:~A \in \mathcal{A}^*(s-2)\right\}. \label{linel}
\end{eqnarray}
Furthermore, the point set $\left\{\mathbf{y}_1^{(s-2)}, \dots, \mathbf{y}_{m-1}^{(s-2)} \right\}$ has property $\mathbf{P}(r)$, which is inherited from the original point set $\{\mathbf{y}_1, \dots, \mathbf{y}_{m}\}$. Thus the term in line \eqref{linel} is dominated by $\hat{N} (m-1,s-2,r)$, and in total we have
\begin{equation} \label{nmsr}
\hat{N} (m,s,r) \leq \hat{N} (m-1,s,r) + \hat{N} (m-1,s-2,r).
\end{equation}
This is an analogue of \eqref{rec}, except that now we ``lose'' two dimensions rather than only one, and that it is only valid as long as $r < s$.\\

Note that from the definition of $\hat{N} (m,s,r)$ we have $\hat{N} (m,s,r) \leq N(m,s)$ for all $m,s,r$. Now we claim the following:\\

{\bf Claim: } We have $\hat{N} (n,d,r) \leq N(n,r) \sum_{0 \leq i \leq d/2} \binom{n}{i}$.\\

The claim is obviously right whenever $r \geq d$, since then 
$$
 \hat{N} (n,d,r) \leq N(n,d) \leq N(n,r).
$$ 
On the other hand, whenever $r < d$, then by \eqref{nmsr} we have
\begin{eqnarray}
\hat{N} (n,d,r) & \leq & \hat{N} (n-1,d,r) + \hat{N} (n-1,d-2,r) \nonumber\\
&\leq &  N(n,r) \sum_{0 \leq i \leq d/2} \binom{n-1}{i} + N(n,r) \sum_{0 \leq i \leq d/2-1} \binom{n-1}{i} \nonumber\\
& = & N(n,r) \sum_{0 \leq i \leq d/2} \left(  \binom{n-1}{i} + \binom{n-1}{i-1} \right) \label{linek} \\
& = & N(n,r) \sum_{0 \leq i \leq d/2} \binom{n}{i}. \label{lineq}
\end{eqnarray}
where in line \eqref{linek} we read $\binom{n-1}{-1}=0$. Thus the claim is true by induction. Classically we have 
$$
\sum_{i=0}^d \binom{n}{i} \leq \left( \frac{en}{d}\right)^d
$$
for $n \geq d$ (see for example \cite[Corollary 3.3]{mohri}), so by \eqref{sauer2} and \eqref{lineq} we have
$$
\hat{N} (n,d,r) \leq  \left( \frac{en}{r}\right)^r  \left( \frac{en}{d/2}\right)^{d/2},
$$
for $n \geq d$, which in particular yields 
\begin{equation} \label{nbound}
\hat{N} (n,d,d/4) \leq  \left( \frac{4en}{d}\right)^{d/4}  \left(\frac{2en}{d}\right)^{d/2} = 2^d \left( \frac{en}{d} \right)^{3d/4}.
\end{equation}

The remaining part of the proof of Theorem \ref{th2} can be carried out similar to the proof of the main theorem in \cite{hinrichs}. As shown in equation (8) of \cite{hinrichs}, for given $\varepsilon>0$ there exists a collection $\mathcal{C}$ of at least $(8 e \varepsilon)^{-d}$ anchored axis-parallel boxes in $[0,1]^d$ such that
$$
\textup{vol} (C_1 \Delta C_2) \geq \varepsilon \qquad \text{for all $C_1,C_2 \in \mathcal{C}$},
$$
where $\Delta$ denotes the symmetric difference. Let $\mathbf{x}_1, \dots, \mathbf{x}_d$ denote the points from the assumption of Theorem \ref{th2}. Since $\mathcal{C}$ is a subset of $\mathcal{A}^*$, by \eqref{nbound} we have
$$
\# \left\{ C \cap \{\mathbf{x}_1, \dots, \mathbf{x}_n\}: ~C \in \mathcal{C} \right\} \leq 2^d \left( \frac{en}{d} \right)^{3d/4}.
$$
Thus by the pigeon hole principle there exist two sets $C_1$ and $C_2$ for which 
\begin{equation} \label{c1c2}
C_1 \cap \{\mathbf{x}_1, \dots, \mathbf{x}_n\} = C_2 \cap \{\mathbf{x}_1, \dots, \mathbf{x}_n\} \qquad \text{and} \qquad \textup{vol} (C_1 \Delta C_2) \geq \varepsilon,
\end{equation}
provided that
\begin{equation} \label{provide}
2^d \left( \frac{en}{d} \right)^{3d/4} < \left(\frac{1}{8 e \varepsilon} \right)^d.
\end{equation}
It is easily seen that \eqref{c1c2} implies that 
$$
D_n^*(\mathbf{x}_1, \dots, \mathbf{x}_n ) \geq \frac{\varepsilon}{4}
$$
(see \cite[Lemma 6]{hinrichs}), and that \eqref{provide} is satisfied if we choose
$$
\varepsilon = \frac{d^{3/4}}{93 n^{3/4}}.
$$
This proves Theorem \ref{th2}.

%%%%%%%%%%%%%%%%%%%%%%%%%%%%%%%%%%%%%%%%%%%%%%%%%%%%%%%%%%%%%%%%%%%%%%%%%%%%%%%%%%%%%%%%%%%
%%% The acknowledgements
\section*{Acknowledgement}

The first author is supported by the Austrian Science Fund (FWF), projects F5507-N26 and I1751-N26, and by the FWF START project Y-901-N35.

%%%%%%%%%%%%%%%%%%%%%%%%%%%%%%%%%%%%%%%%%%%%%%%%%%%%%%%%%%%%%%%%%%%%%%%%%%%%%%%%%%%%%%%%%%%
%%% The bibliography
%
% BibTeX users please use
%\bibliographystyle{spmpsci}
%\bibliography{mybiblio}
%   and then copy and paste the content of the .bbl file here for the final version.
%
% E.g.:

%\bibliography{Cheap_proof}
%\bibliographystyle{spmpsci}

\end{document}